\newif\ifsubmit
\submitfalse

\ifsubmit
    \documentclass[letterpaper, 10 pt, conference]{ieeeconf}  

    \IEEEoverridecommandlockouts                              
\else
    \documentclass[12pt, hidelinks]{article}
\fi
\usepackage[utf8]{inputenc}
\usepackage{amsmath, amssymb, url, float, graphicx, float, xcolor}
\usepackage{fullpage}
\usepackage[small,bf]{caption}
\usepackage{adjustbox}
\usepackage{hyperref}
\usepackage{subcaption}

\ifsubmit\else
\usepackage{jlcode}
\def\jlbasicfont{\ttfamily\footnotesize\selectfont}
\lstdefinestyle{code}{
  language=Julia, 
  showstringspaces=false,
  keywordstyle=\color{blue},
  commentstyle=\color{gray},
  identifierstyle=\color[RGB]{0,102,0},
  columns=fullflexible,
  keepspaces=true
}
\lstnewenvironment{code}{\lstset{style=code}}{}

\newcommand{\codelink}[1]{\lstinputlisting[style=code]{#1}
  \noindent\begin{center}
  \filename@parse{#1}
  \href{\codeurl/\filename@base.\filename@ext}
  {\textcolor{blue}{\codelinktext}}
  \end{center}
}
\newcommand{\linkonly}[1]{
  \protect\filename@parse{#1}
  \href{\codeurl/\filename@base.\filename@ext}
  {\textcolor{blue}{\tt[\protect\filename@base.\protect\filename@ext]}}
}

\fi

\usepackage[style=alphabetic,backend=bibtex]{biblatex}
\bibliography{refs}

\usepackage{tikz, pgfplots}
\pgfplotsset{
    compat=1.16,
}
\usetikzlibrary{topaths,intersections,calc}
\usetikzlibrary{positioning,chains,fit,shapes,calc,arrows.meta}

\let\phi\varphi
\let\eps\varepsilon

\newcommand{\BEAS}{\begin{eqnarray*}}
    \newcommand{\EEAS}{\end{eqnarray*}}
    \newcommand{\BEA}{\begin{eqnarray}}
    \newcommand{\EEA}{\end{eqnarray}}
    \newcommand{\BEQ}{\begin{equation}}
    \newcommand{\EEQ}{\end{equation}}
    \newcommand{\BIT}{\begin{itemize}}
    \newcommand{\EIT}{\end{itemize}}
    \newcommand{\BNUM}{\begin{enumerate}}
    \newcommand{\ENUM}{\end{enumerate}}
    
    \newcommand{\cf}{{\it cf.}}
    \newcommand{\eg}{{\it e.g.}}
    \newcommand{\ie}{{\it i.e.}}

    \newcommand{\ones}{\mathbf 1}
    
    \newcommand{\reals}{\mathbf{R}}

    \newcommand{\bmat}[1]{\begin{bmatrix}#1\end{bmatrix}}

    \newif\iftodos

\newcommand{\solver}{\texttt{ConvexFlows}}

\ifsubmit
    \title{\Large \bf Solving the Convex Flow Problem}
    \author{Theo Diamandis and Guillermo Angeris
    \thanks{T.\ Diamandis is supported by the Department of Defense (DoD)
    through the National Defense Science \& Engineering Graduate (NDSEG)
    Fellowship Program.}
    \thanks{T.\ Diamandis is with Computer Science and Artificial Intelligence Lab,
            Massachusetts Insttitue of Technology, Cambrigde, MA, USA
            {\tt\small tdiamand@mit.edu}}%
    \thanks{G.\ Angeris is with Bain Capital,
            Boston, MA, USA
            {\tt\small gangeris@baincapital.com}}%
}
\else
    \title{Solving the Convex Flow Problem}
    \author{Theo Diamandis \and Guillermo Angeris}
    \date{March 2024}
\fi

\begin{document} 
\maketitle
\ifsubmit
    \thispagestyle{empty}
    \pagestyle{empty}
\fi

\begin{abstract}
    In this paper, we introduce the solver \solver{} for the convex flow problem
    first defined in the authors' previous work. In this problem, we aim to
    optimize a concave utility function depending on the flows over a graph.
    However, unlike the classic network flows literature, we also allow for a
    concave relationship between the input and output flows of edges. This
    nonlinear gain describes many physical phenomena, including losses in power
    network transmission lines. We outline an efficient algorithm for solving
    this problem which parallelizes over the graph edges. We provide an open
    source implementation of this algorithm in the Julia programming language
    package \texttt{ConvexFlows.jl}. This package includes an interface to
    easily specify these flow problems. We conclude by walking through an
    example of solving for an optimal power flow using \solver{}.
\end{abstract}

\section{Introduction}
Theorists and practitioners both apply network flow models to describe, analyze,
and solve problems from many domains---from routing trucks to routing bits. For
linear flows, an extensive academic literature developed the associated theory,
algorithms, and applications. (See, \eg,~\cite{ahuja1988network},
\cite{williamson2019network}, and references therein.) However, these linear
models often fail to describe real systems. For example, in electrical systems,
the power lost increases as more power is transmitted; in communications
systems, the message failure rate increases as more messages are transmitted;
and, in financial systems, the price of an asset increases as more of that
asset is purchased. In each of these cases, the output of the system is a
concave function of its input.

In this work, we focus on solving this more general \emph{convex flow problem},
an important special case of the authors' previous
work~\cite{diamandis2024convexflows}, and provide a package with a clean
interface to do so. Although this problem is a convex optimization problem, for
which many open-source and commercial solvers exist, the convex flow problem
has additional structure that can be exploited. Following this prior
work~\cite{diamandis2024convexflows}, we use a dual decomposition approach,
which allows us to decompose the problem over the network edges. In contrast
with the previous approach, though, we solve this problem using the
Broyden–Fletcher–Goldfarb–Shanno (BFGS)
method~\cite[\S6]{nocedalNumericalOptimization2006}. This method has been shown
to be robust against non-smooth objective functions~\cite{lewis2013nonsmooth}
that often appear in practical instances of the convex flow problem.  To specify
these problems, we provide an easy-to-use interface that, unlike in previous
work, does not require specifying conjugate functions or the support functions
for feasible sets. We provide an open-source implementation of the method and
this interface in the Julia programming language\ifsubmit\footnote{ Available online at
\url{https://github.com/tjdiamandis/ConvexFlows.jl}.}\else\fi with extensive
documentation. We conclude with two optimal power flow examples and associated
numerical experiments. Additional examples are available in the \solver{}
documentation. \ifsubmit\else All code is available online at 
\begin{center}
    \texttt{\url{https://github.com/tjdiamandis/ConvexFlows.jl}.}
\end{center}
\fi

\section{The convex flow problem}\label{sec:convex-flow}
We consider a directed graph with $n$ nodes and $m$ edges. Each edge, $i = 1, \dots, m$, in
the graph has an associated strictly concave,
nondecreasing gain function $h_i: \reals_+ \to \reals_+ \cup \{-\infty\}$, which
denotes the output flow $h_i(z)$ of edge $i$ given some input flow $z \in
\reals_+$. (We assume strict concavity, but this can be achieved generally by,
say, subtracting a small quadratic term from the gain function.) We use
infinite values to encode constraints: an input flow $z$ over edge $i$ such
that $h_i(z) = -\infty$ is unacceptable. We denote the \emph{flow} across the
edge by the vector $x_i \in \reals^2$, where $x_1 \le 0$ is the flow into edge
(equivalently, out of edge $i$'s source node) $i$ and $x_2 \ge 0$ is the flow
out of the edge (equivalently, into edge $i$'s terminal node). These flows are
connected by the relationship
\[
    x_2 = h_i(-x_1).
\]
With each edge $i$ we associate a matrix $A_i \in \{0,1\}^{n \times 2}$ that
maps the `local' indices of nodes to their global indices. More specifically,
if edge $i$ connects node $j$ to node $j'$, then we define $A_i = \bmat{e_j &
e_{j'}}$, where $e_j$ denotes the $j$th unit basis vector. After mapping each
edge flow to the global index and summing across all edges, we obtain the
\emph{net flow vector} $y \in \reals^n$, defined as
\[
    y = \sum_{i=1}^m A_i x_i.
\]
If $y_j > 0$, then this node has flow coming into it and is called a
\emph{sink}. If $y_j < 0$, then this node provides flow to the network and is
called a \emph{source}.

In the convex flow problem, we aim to maximize some utility function $U:
\reals^n \to \reals \cup \{-\infty\}$ over all feasible net flows $y$. Infinite
values again denote constraints: any flow with $U(y) = -\infty$ is unacceptable.
We require this utility function to be strictly concave and strictly increasing.
(The nondecreasing utility case also follows directly from this setup but
requires some additional care.) The convex flow problem is
\begin{equation}\label{eq:primal-problem}
    \begin{aligned}
        & \text{maximize} && U(y)\\
        & \text{subject to} && {\textstyle y = \sum_{i=1}^m A_ix_i}\\
        &&& (x_i)_2 \le h_i\left(-(x_i)_1\right), ~~ i=1, \dots, m.
    \end{aligned}
\end{equation}
An important consequence of this setup is that a solution $\{x_i^\star\}$
to~\eqref{eq:primal-problem} will always saturate edge inequality constraints;
\ie, $(x_i^\star)_2 = h(-(x_i^\star)_1)$. To see this, note that any flow $x_i$
satisfying $(x_i)_2 < h\left(-(x_i)_1\right)$, can have its second component
increased to $(x_i)_2 + \eps $ for some $\eps > 0$. Since $U$ is strictly
increasing in $y$, and $y$ is (elementwise) strictly increasing in the $x_i$,
this change would increase the objective value, so these flows $x_i$ could not
have been optimal.

\begin{figure}
    \centering
    \hfill
    \adjustbox{max width=0.48\columnwidth}{
        \input figures/opf-gain.tex
    }
    \hfill
    \adjustbox{max width=0.48\columnwidth}{
        \begin{tikzpicture}[scale=1.0]
    \begin{axis}[
        axis on top=true,
        xmin=-4.5, xmax=0.5,
        ymin=-0.1, ymax=4.5,
        axis lines=center,
        xlabel={$z_1$},
        x label style={anchor=north},
        ylabel={$z_2$},
        grid=major,
        legend pos=outer north east,
        samples=200,
        domain=-4:4,
        xtick=\empty,
        ytick=\empty
        ]
        \addplot[black, very thick, domain=-4:0] {3*-x - 16*(ln(1 + exp(-x/4)) - ln(2))};
        \addplot[gray!30, domain=-4:0, fill, opacity=0.6, draw=none] {3*-x - 16*(ln(1 + exp(-x/4)) - ln(2))} -| (-4, 0) -- cycle;
        \node at (-2.5, 0.75) {$T$};

        \draw[black, very thick] (-4, 2.09) -- (-4, 0);
        \draw[black, very thick] (0, 0) -- (-4, 0);
    \end{axis}
\end{tikzpicture}
    }
    \hfill \null
    \caption{
        A concave gain function $h$ with implicitly bounded domain (left), and
        the corresponding set of allowable flows (right).
    }
    \label{fig:opf-loss}
\end{figure}
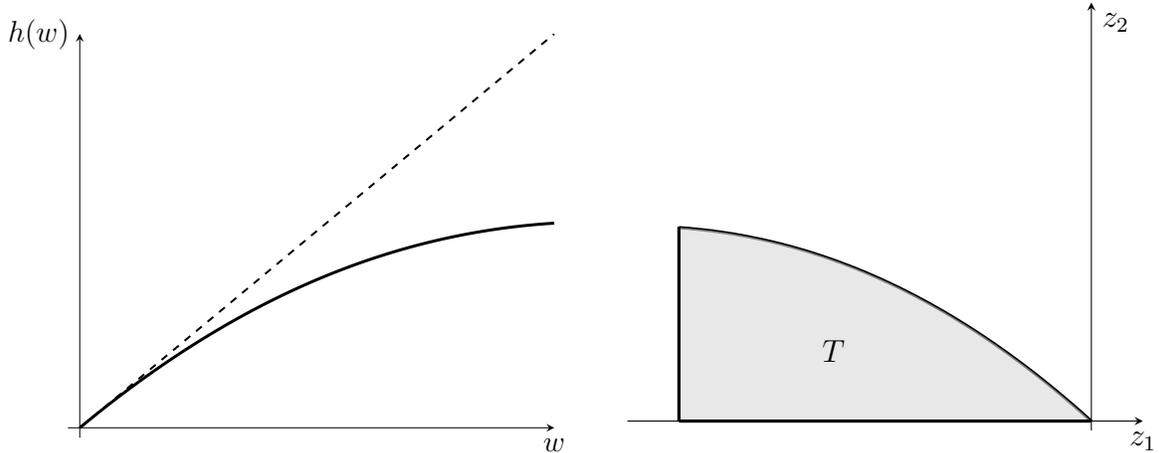

In what follows, we will call a flow $x_i$ over edge $i$ an \emph{allowable
flow} if it satisfies the constraint in~\eqref{eq:primal-problem}:
\[
    (x_i)_2 \le h_i(-(x_i)_1).
\]
Note that in prior work~\cite{diamandis2024convexflows}, we instead defined a
closed convex set of allowable flows $T_i$ for each edge $i$. This set can be
constructed directly from the inequality above; figure~\ref{fig:opf-loss} shows
an example.

\section{Dual problem}
Observe that the convex flow problem~\eqref{eq:primal-problem} only has one
constraint coupling the edge flows $x_i$. This structure suggests that we should
relax the linear equality constraint to a penalty and consider the resulting
dual problem~\cite[\S5.2]{cvxbook}:
\begin{equation}
    \label{eq:dual-problem}
    \begin{aligned}
        & \text{minimize} && \bar U(\nu)
            + \sum_{i=1}^m f_i(A_i^T \nu).
    \end{aligned}
\end{equation}
The only variable in this problem is $\nu \in \reals^n$ and the
functions $\bar U$ and $f_i$ are defined
\begin{subequations}
    \label{eq:subproblems}
    \begin{align}
        \label{eq:subproblem-U}
        \bar U(\nu) &= \sup_{y} (U(y) - \nu^Ty), \\
        \label{eq:subproblem-f}
        f_i(\eta) &= \sup_{w \ge 0} ( -\eta_1w + \eta_2h(w) ),
    \end{align}
\end{subequations}
for $i=1, \dots, m$. Note that $\bar U(\nu) = (-U)^*(-\nu)$ is the is the
Fenchel conjugate~\cite[\S3.3]{cvxbook} of $-U$ with a negated argument,
while $f_i$ is essentially the support function for the set
\[
    \{(z_1, z_2) \mid z_2 \le h_i(-z_1)\},
\]
if $\eta \ge 0$. This fact follows from problem~\eqref{eq:dual-problem} if
$A_i^T\nu \ge 0$ for all $i$, or, equivalently, if $\nu \ge 0$, which we show
next.

\subsection{Properties}
Assuming that there exists a point in the relative interior of the feasible set
(Slater's condition), the dual problem~\eqref{eq:dual-problem} has the same
optimal value as the primal problem~\eqref{eq:primal-problem}. This assumption
typically holds in practice, so we will focus on
solving~\eqref{eq:dual-problem}. We will show two things: first, that any
optimal $\nu^\star$ is nonnegative (and, indeed, that $\nu^\star > 0$) since $U$ is
strictly increasing, and, second, that, given $\nu$
solving~\eqref{eq:dual-problem}, the solutions to the
subproblems~\eqref{eq:subproblems} are feasible for the primal problem and
therefore optimal. The first fact will be useful in solving the dual problem,
while the second fact will imply that, by solving the dual
problem~\eqref{eq:dual-problem}, we can recover a solution to the original
problem~\eqref{eq:primal-problem}.

First, let $y$ be any point with $U(y) > -\infty$. If $\nu_j < 0$ for some $j$,
we would have
\ifsubmit
\begin{multline*}
    \bar U(\nu) \ge U(y+te_j) - (y + te_j)^T\nu \\
    \ge U(y) - \nu^Ty - t\nu_j \to \infty.
\end{multline*}
\else
\[
    \bar U(\nu) \ge U(y+te_j) - (y + te_j)^T\nu 
    \ge U(y) - \nu^Ty - t\nu_j \to \infty.
\]
\fi
as $t \to \infty$, where, in the second inequality, we have used the fact that
$U$ is increasing. Therefore, for $\bar U(\nu)$ to be finite, we must have $\nu
\ge 0$. We will show soon that the second claim implies that any optimal dual
variables satisfy $\nu > 0$, if the primal problem~\eqref{eq:primal-problem}
has a finite solution.

For the second claim: it is not hard to show that $\bar U$ and $f_i$ are
differentiable, when finite, since $U$ and the $h_i$ are strictly
concave~\cite[Thm.\ 25.1]{rockafellar1970convex}. Let $\nu^\star$
be dual optimal, then the first order optimality
conditions applied to problem~\eqref{eq:dual-problem} give
\begin{equation}
    \label{eq:dual-grad}
    -y^\star(\nu^\star) + \sum_{i=1}^m A_ix_i^\star(\nu^\star) = 0,
\end{equation}
where 
\[
    y^\star(\nu) = \nabla \bar U(\nu)
\] 
is the maximizer for
subproblem~\eqref{eq:subproblem-U} and 
\[
    x_i^\star(\nu) = \nabla f_i(A_i^T\nu) = (-w^\star(\nu), h(w^\star(\nu))),
\] 
where $w^\star$ is the maximizer for subproblem~\eqref{eq:subproblem-f}. Since
these points are feasible for~\eqref{eq:primal-problem} then they must also be
optimal.

Finally, if~\eqref{eq:primal-problem} has a finite solution $y^\star$, then the
first-order optimality condition for~\eqref{eq:subproblem-U} means that
$\nabla U(y^\star) = \nu^\star$. But, since $U$ is strictly increasing, we have
that $\nabla U(y^\star) > 0$ so $\nu^\star > 0$ as required.

\subsection{Solving the dual problem}
The fact that $\nu^\star > 0$ suggests a natural way of modifying a solution
method to respect this constraint: we simply modify a line search
to ensure $\nu$ remains positive. Specifically, we add an upper bound on the
step size which ensures that every iterate remains strictly positive. This approach
keeps the problem otherwise unconstrained, which simplifies solution methods.

For small to medium-sized problems, we use the quasi-Newton method BFGS, which
has been shown to work well for nonsmooth problems~\cite{lewis2013nonsmooth}. We
use the bracketing line search from Lewis and Overton~\cite{lewis2013nonsmooth},
modified to prevent steps outside of the positive orthant, which also ensures
that the step size satisfies the weak Wolfe conditions. (We note that \solver{}
also includes an interface to
L-BFGS-B~\cite{byrd1995limited,zhu1997algorithm,morales2011remark} for larger
problems, but this interface requires more a more sophisticated problem
specification, and this method may be less robust to nonsmoothness in the 
problem~\cite{asl2021behavior}.)

Importantly, evaluating the dual objective function and its
gradient~\eqref{eq:dual-grad} parallelizes across all the edges, and each
individual subproblem can be solved very quickly---often in closed form. This
observation suggests a natural interface to specify the convex flow problem: we
only need a means of evaluating the subproblems and computing their maximizers.
Given user-specified utility and gain functions, our software automatically
computes these subproblem evaluations.

\section{Interface}
It is unreasonable to expect most users to directly specify conjugate-like
functions and solutions to convex optimization problems as
in~\eqref{eq:subproblem-U} and~\eqref{eq:subproblem-f}. Instead, we develop an
interface that allows the user to specify the utility function $U$ and the edge
gain functions $h_i$ for each edge $i = 1, \dots, m$. With this, and the
previous discussion, we can now solve the dual problem and, from there, recover
a primal optimal solution.

\subsection{The first subproblem}
The first subproblem~\eqref{eq:subproblem-U} typically has a closed form
expression. Since, from before, $\bar U(\nu) = (-U)^*(-\nu)$, where $U^*$
denotes the Fenchel conjugate of $U$, we can use standard results in conjugate
function calculus to compute $\bar U$ from a number of simpler `atoms'. For
example, $U$ is often separable, in which case we have that $U(y) = u_1(y_1) +
\dots + u_n(y_n)$, so
\[
    \bar U(y) = \bar u_1(y_1) + \dots + \bar u_n(y_n),
\]
where $\bar u_j$ is defined identically to~\eqref{eq:subproblem-U}. Our package
\solver{} provides atoms that a user can use to construct $U$. Some examples of
scalar utility atoms include the linear, nonnegative linear, and nonpositive
quadratic atoms. We also provide a number of cost functions, including
nonnegative quadratic cost. Note that, since $U$ is increasing, we can support
lower bounds on the variables but not upper bounds.

While it is most efficient to build $U$ (and therefore $\bar U$) from known
atoms, more general functions without constraints may be handled by
solving~\eqref{eq:subproblem-U} directly. A vector $\tilde y$ achieving the
supremum must satisfy $\nabla U(\tilde y) = \nu$. This equation may be solved
via Newton's method, and the gradient and Hessian may be computed via
automatic differentiation.

We can also incorporate constraints by writing $U$ as the solution to a conic
optimization problem, which may be expressed using a modeling language such as
JuMP~\cite{dunningJuMPModelingLanguage2017,Lubin2023} or
\texttt{Convex.jl}~\cite{udell2014convex}, both of which can compile problems
into a standard conic form using
\texttt{MathOptInterface.jl}~\cite{legatMathOptInterfaceDataStructure2021}.

\subsection{The second subproblem}
For each edge $i$ we require the user to specify the gain function $h_i$ in
native Julia code. Denote the solution point of the second
problem~\eqref{eq:subproblem-f} by $w^\star$. We write $h^+(w)$ and $h^-(w)$
for the right and left derivatives of $h$ at $w$, respectively. Specifically,
we define
\[
    h^+(w) = \lim_{\delta \to 0^+} \frac{h(w + \delta) - h(w)}{\delta},
\]
and $h^-(w)$ analogously. The
optimality conditions for problem~\eqref{eq:subproblem-f} are then that
$w^\star$ is a solution if, and only if,
\begin{equation}\label{eq:swap-opt-cond}
     h^+(w^\star) \le \eta_1/\eta_2 \le  h^-(w^\star).
\end{equation}
(We may assume $\eta_2 > 0$ from the previous discussion, since $\nu > 0$.)
Note that the optimality condition suggests a simple method to check if an edge
will be used at all: zero flow is optimal if and only if
\[
    h^+(0) \le {\eta_1}/{\eta_2} \le h^-(0).
\]
This `no flow condition' is often much easier to check in practice than solving
the complete subproblem.

If the zero flow is not optimal, then we can solve~\eqref{eq:subproblem-f} via a
one-dimensional root-finding method. We assume that $h$ is differentiable almost
everywhere (\eg, $h$ is a piecewise smooth function) and use bisection search or
Newton's method to find a $w^\star$ that satisfies~\eqref{eq:swap-opt-cond}.
Since we use directed edges, and typically an upper bound $b$ on the flow exists
for physical systems, we begin with the bounds $(0, b)$ and terminate after
$\log_2(b/\eps)$ iterations. (If no bound is specified, an upper bound $b$ may
be computed with, for example, a doubling method.) We compute the first
derivative of $h$ using forward mode automatic differentiation, implemented in
\texttt{ForwardDiff.jl}~\cite{RevelsLubinPapamarkou2016}. Computing a
derivative can be done simultaneously with a function evaluation and, as a
result, these subproblems can typically be solved very quickly. Alternatively,
the user may specify a closed-form solution to the subproblem, which exists for
many problems in practice (see, for example, the examples
in~\cite[\S6]{diamandis2024convexflows}.)

\section{Example: optimal power flow}
We adapt the optimal power flow example
of~\cite[\S3.2]{diamandis2024convexflows}. This problem seeks to find a
cost-minimizing plan to generate power, which may be transmitted over a network
of $m$ transmission lines, to satisfy the power demand of $n$ regions over some
number of time periods $T$. We use the transport model for power networks along
with a nonlinear transmission line loss function from~\cite{stursberg2019mathematics},
which results in a good approximation of the DC power flow model.

The loss function models the phenomenon that, as more power is transmitted along
a line, the line dissipates an increasing fraction of the power transmitted.
Following~\cite[\S2]{stursberg2019mathematics}, we use the convex, increasing
loss function
\[
    \ell_i(w) = \alpha_i \left(\log(1 + \exp(\beta_i w)) - \log 2\right) - 2w,
\]
where $\alpha_i$ and $\beta_i$ are known constants for each line and satisfy
$\alpha_i\beta_i = 4$. The gain function of a line with input $w$ can then be
written as
\[
    h_i(w) = w - \ell_i(w).
\] 
Each line $i$ also has a maximum capacity, given by $b_i$.
Figure~\ref{fig:opf-loss} shows a power line gain function and its corresponding
set of allowable flows.

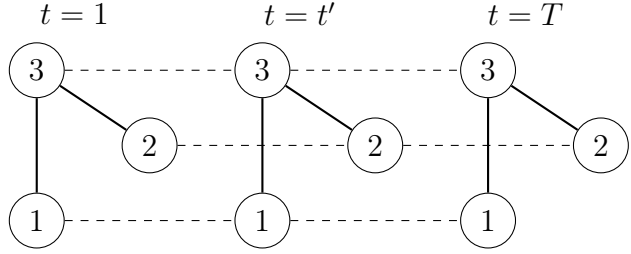
\begin{figure}[t]
    \centering
    \hfill
    \adjustbox{max width=0.95\columnwidth}{
        \begin{tikzpicture}
\node[circle, draw] (1) at (0,0) {3};
\node[circle, draw] (2) at (1.5,-1.0) {2};
\node[circle, draw] (3) at (0,-2) {1};

\draw[thick] (1) -- (2);
\draw[thick] (1) -- (3);

\node[circle, draw] (1t) at (3,0) {3};
\node[circle, draw] (2t) at (4.5,-1.0) {2};
\node[circle, draw] (3t) at (3,-2) {1};

\draw[thick] (1t) -- (2t);
\draw[thick] (1t) -- (3t);

\node[circle, draw] (1T) at (6,0) {3};
\node[circle, draw] (2T) at (7.5,-1.0) {2};
\node[circle, draw] (3T) at (6,-2) {1};

\draw[thick] (1T) -- (2T);
\draw[thick] (1T) -- (3T);

\draw[dashed] (1) -- (1t);
\draw[dashed] (2) -- (2t);
\draw[dashed] (3) -- (3t);

\draw[dashed] (1t) -- (1T);
\draw[dashed] (2t) -- (2T);
\draw[dashed] (3t) -- (3T);

\node at (0.5,0.75) {$t = 1$};
\node at (3.5,0.75) {$t = t'$};
\node at (6.5,0.75) {$t = T$};

\end{tikzpicture}
    }
    \hfill \null
    \caption{
        Graph representation of a power network with three nodes over time. 
        Each solid line corresponds to a transmission line edge, and each dashed
        line corresponds to a storage edge.
    }
    \label{fig:opf-graph}
\end{figure}

Each node $j$ may also store power generated at time $t$ for use at time
$t+1$. If $w$ units are stored, then $\gamma_j w$ units are available at time
$t+1$ for some $\gamma_j \in [0, 1]$. These parameters may describe, for example,
the battery storage efficiency. We model this setup by introducing $T$ nodes in the graph
for each node, with an edge from the $t$th node to the $(t+1)$th node
corresponding to node $j$ with the appropriate linear gain function, as
depicted in figure~\ref{fig:opf-graph}. (Note that, for numerical stability, we
subtract a small quadratic term, $(\eps/2)w^2$, from the linear gain functions,
where $\eps$ is very small.)

At time $t = 1, \dots, T$, node $j = 1, \dots, n$ demands $d_{tj}$ units of power and can generate
power $p_j$ at a cost $c_j : \reals \to \reals_+$, given by
\[
    c_j(p) = \begin{cases}
        (\kappa_j/2)p^2 & p \ge 0 \\
        0 & p < 0,
    \end{cases}
\]
which is a convex, increasing function parameterized by $\kappa_j > 0$. Power
dissipation has no cost but also generates no profit. To meet demand, we must
have that, for each $t = 1, \dots, T$,
\[
    d_t = p_t + y_t, \qquad \text{where} \qquad  y_t =  \sum_{i=1}^m A_i x_{ti}.
\]
In other words, the power produced, plus the net flow of power, must satisfy the 
demand in each node.
We write the network utility function as
\begin{equation}\label{eq:app-opf-U}
    U(y) = \sum_{t=1}^T\sum_{j=1}^n -c_j(d_{tj} - y_{tj}).
\end{equation}
Since $c_i$ is convex and nondecreasing in its argument, the utility function
$U$ is concave and nondecreasing in $y$. This problem can then be cast as a
special case of~\eqref{eq:primal-problem}.

Note that the subproblems associated with the optimal power flow problem may be
worked out in closed form. The first subproblem is
\[
    \bar U(\nu) = \sum_{t=1}^T\sum_{j=1}^n\left(\frac{\nu_{tj}^2}{2\kappa_j} - d_{tj} \nu_{tj}\right),
\]
with domain $\nu \ge 0$. The second subproblem is
\[
    f_i(\eta_i) = \sup_{0 \le w \le b_i} \left\{
         -\eta_1 w  + \eta_2 \left( w - \ell_i(w) \right)
    \right\}.
\]
Using the first order optimality conditions, we can compute the solution:
\[
    w^\star_i = \left(\beta_i^{-1} \log\left(\frac{3\eta_2 - \eta_1}{\eta_2 + \eta_1}\right)\right)_{[0, b_i]}, \qquad
\]
where $(\cdot)_{[0, b_i]}$ denotes the projection onto the interval $[0, b_i]$.
These closed form solutions can be directly specified by the user in \solver{}
for increased efficiency.

\subsection{Numerical examples}

\subsubsection{Multi-period power generation example}
We first consider an example network with three nodes over a time period of $5$
days. The first two nodes are users who consume power and have a sinusoidal
demand with a period of $1$ day. These users may generate power at a very high
cost ($\kappa_j = 100$). The third node is a generator, which may generate power
at a low cost ($\kappa_j = 1$) and demands no power for itself. We equip the
second user with a battery, which can store power between time periods with
efficiency $\gamma = 1.0$. For each transmission line, we set $\alpha_i = 16$
and $\beta_i = 1/4$. The network has a total of $360$ nodes and $359$ edges.
\ifsubmit\else The full code can be found in
appendix~\ref{app:multi-period-code}.\fi

\begin{figure}[h!]
    \centering
    \includegraphics[width=\ifsubmit 0.95\else 0.85\fi\columnwidth]{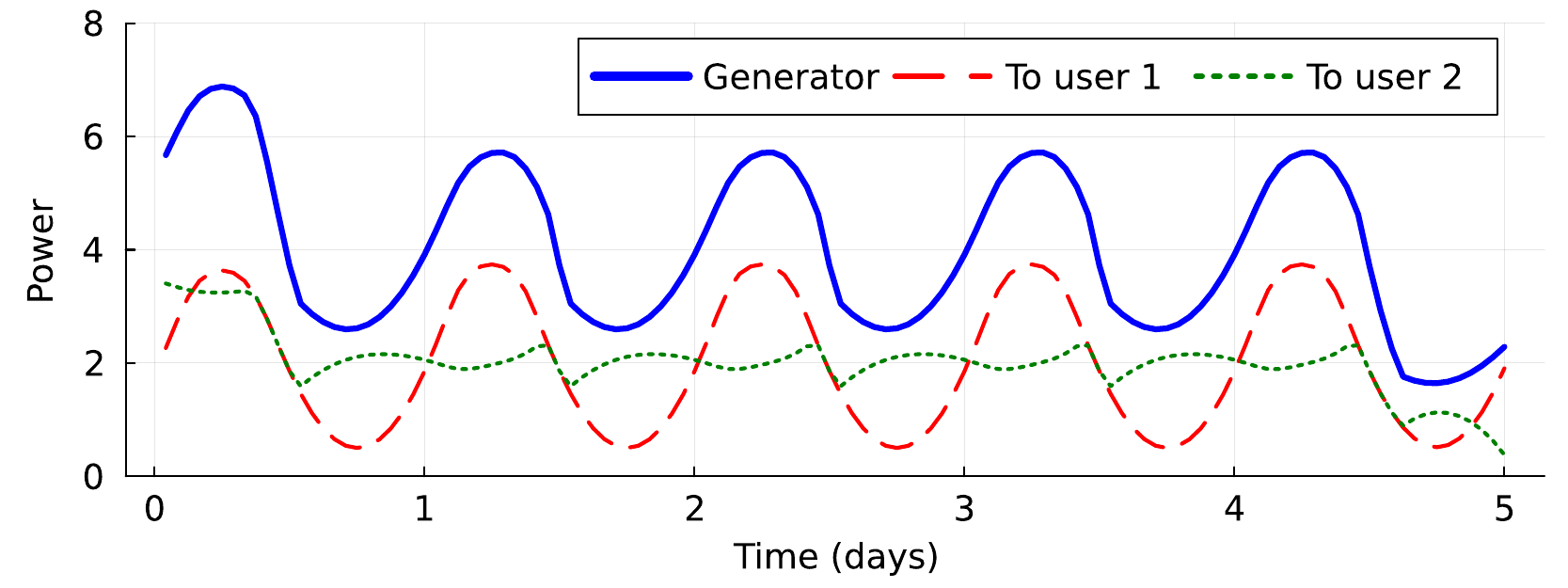}
    \includegraphics[width=\ifsubmit 0.95\else 0.85\fi\columnwidth]{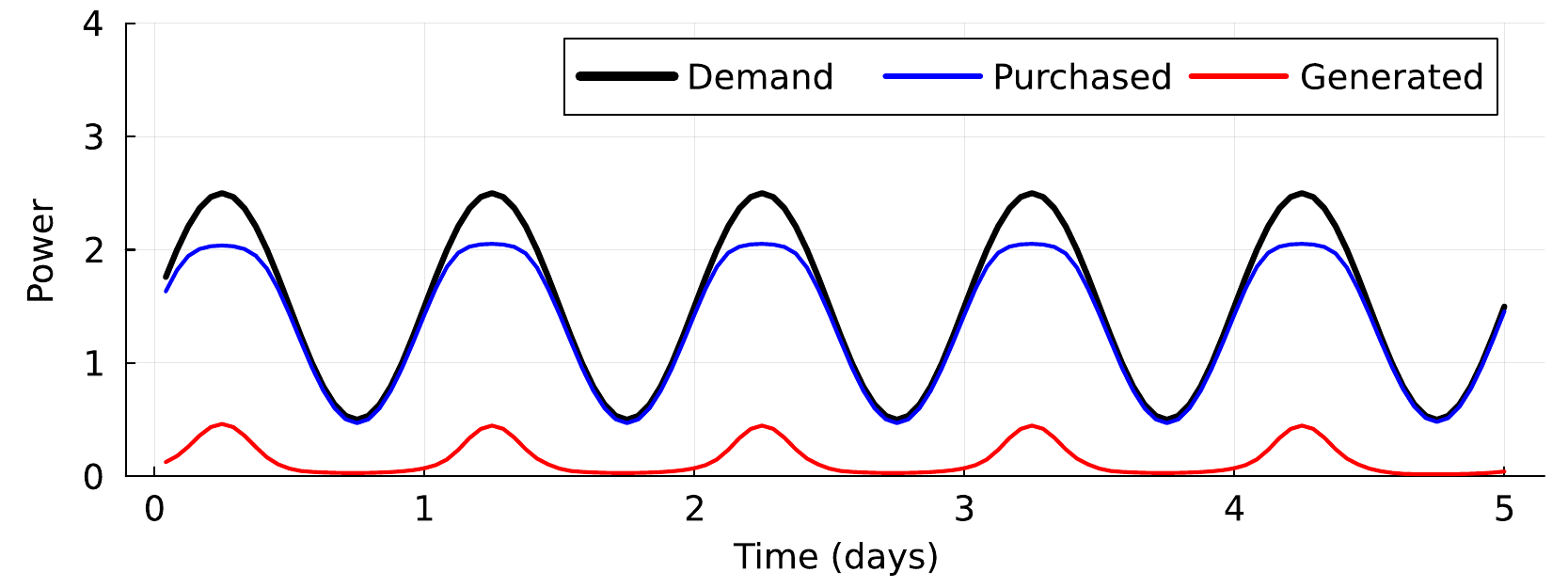}
    \includegraphics[width=\ifsubmit 0.95\else 0.85\fi\columnwidth]{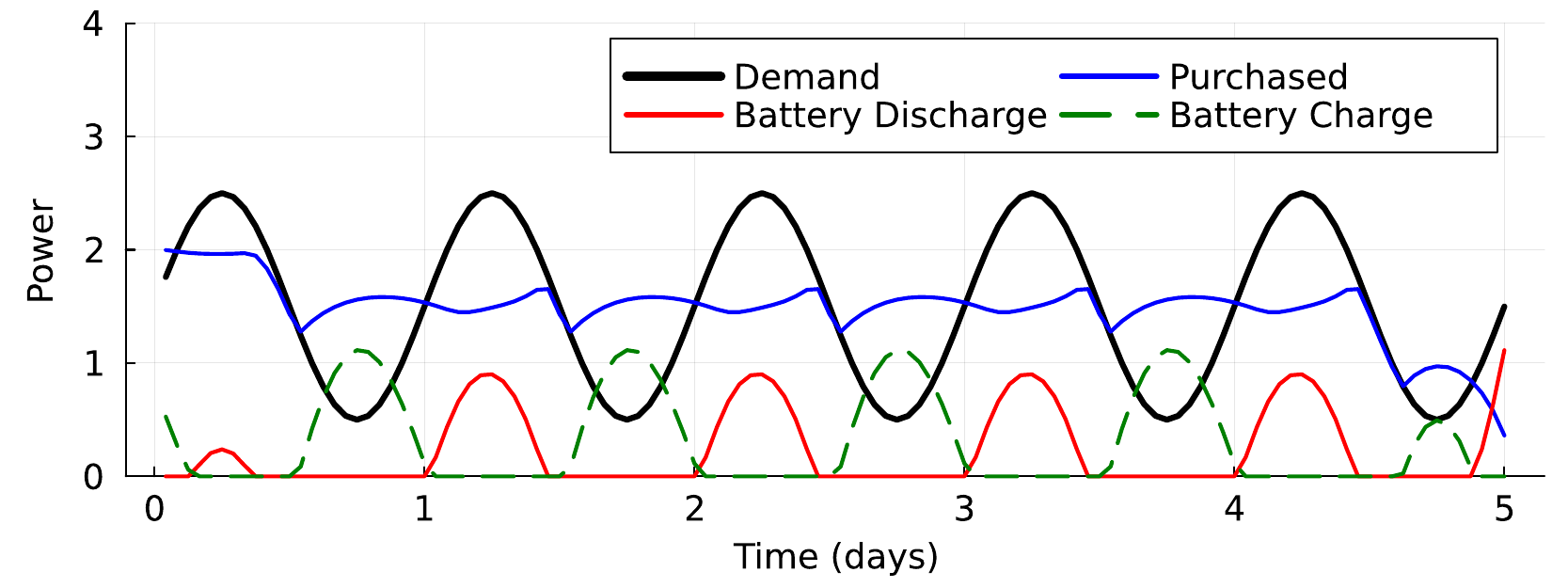}
    \caption{
        Power generated (top), power used by the first node (middle) and by the
        second node, which has a battery (bottom).
    }
    \label{fig:opf-five-day}
\end{figure}

We display the minimum cost power generation schedule in
figure~\ref{fig:opf-five-day}. Notice that during period of high demand, the
first user must generate power at a high cost. The second user, on the other
hand, purchases more power during periods of low demand to charge their battery
and then uses this stored power during periods of high demand. As a result, the
power purchased by this user stays roughly constant over time, after some
initial charging.

\subsubsection{Larger network}
We next consider the network from~\cite{kraning2013dynamic}, generated using the
same parameters. We use $n=100$ nodes and $T = 2$ time periods. For each time
period $t$, we draw the demand $d_{it}$ for each node uniformly at random from
$[1, 5]$. For each transmission line, we again set $\alpha_i = 16$ and $\beta_i
= 1/4$. Each transmission line has maximum capacity drawn uniformly at random
from the set $\{1, 2, 3\}$. A line with maximum capacity $1$ operating at full
capacity will loose about $10\%$ of the power transmitted, whereas a line with
maximum capacity $3$ will loose almost $40\%$ of the power transmitted (\cf,
figure~\ref{fig:opf-loss}). We let all lines be bidirectional: if there is a
line connecting node $j$ to node $j'$, we add a line connecting node $j'$ to
node $j$ with the same parameters. For each node, we allow it to store power
with probability $1/2$ and then draw its efficiency parameter $\gamma_j$
uniformly at random from $[0.5, 1.0]$. In this setup, there are a total of $452$
edges.
\ifsubmit\else See appendix~\ref{sec:interface-opf} for example code.\fi

\begin{figure}[h!]
    \centering
    \includegraphics[width=\ifsubmit 0.95\else 0.85\fi\columnwidth]{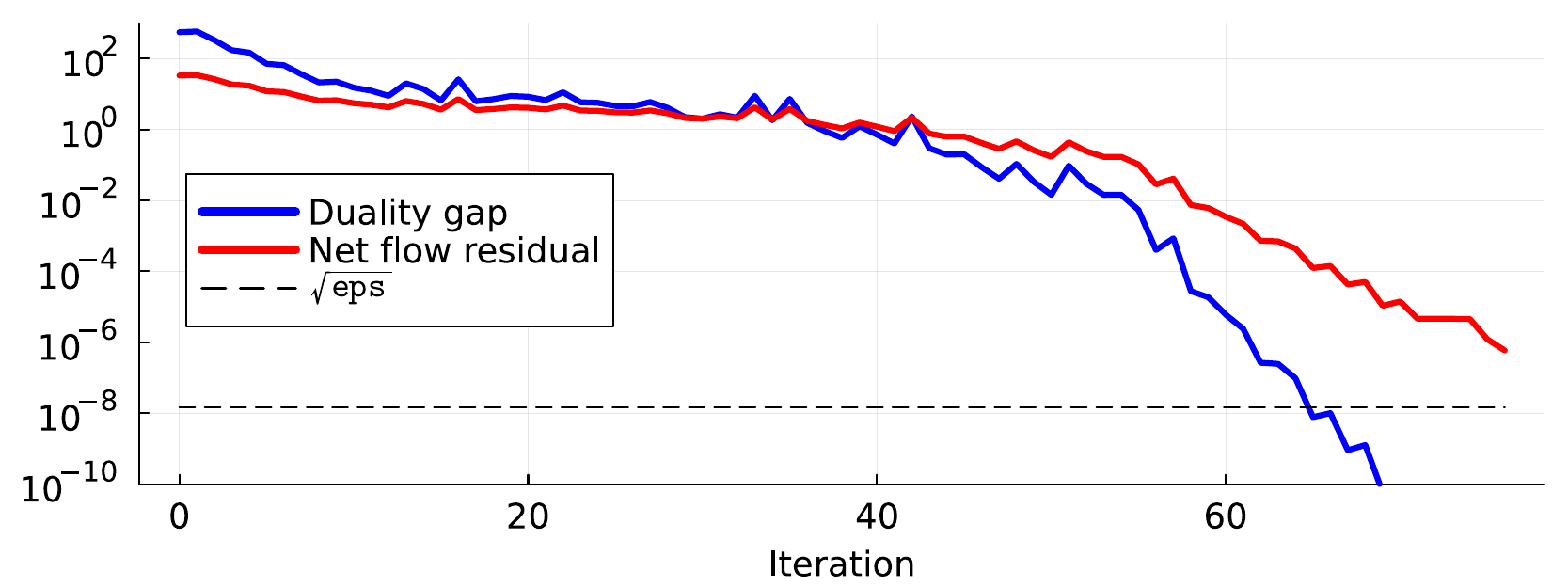}
    \caption{
        Convergence of \texttt{ConvexFlows} with $n = 100$. The primal
        residual measures the net flow constraint violation, with $\{x_i\}$
        from~\eqref{eq:subproblem-f} and $y$ from~\eqref{eq:subproblem-U}.
    }
    \label{fig:opf-two-node-iter}
\end{figure}

Figure~\ref{fig:opf-two-node-iter} shows the convergence of our method on the
optimal power flow problem for this network. The primal feasible point used to
compute the relative duality gap is constructed as
\[
    \hat y = \sum_{i=1}^m A_i \tilde x_i,
\]
where $\tilde x_i$ solves the subproblem~\eqref{eq:subproblem-f} with the
current iterate $\nu_k$. There is a clear linear convergence region, followed by
quadratic convergence, similar to Newton's method. We note that L-BFGS does not
exhibit good convergence on our problem, which is consistent to the results
in~\cite{asl2021behavior}. (See~\cite[\S6.1]{diamandis2024convexflows} for
additional examples using L-BFGS-B to solve the convex flow problem on very 
large networks.)

\section{Conclusion}
This paper introduces the software package \solver{} for solving the convex
flow problem defined in~\cite{diamandis2024convexflows}. This package provides
an easy-to-use interface for specifying these problems, which appear in many
applications, including the transport model optimal power flow problem discussed here. We posit
that additional structure of this problem may be exploited in solution methods.
For example, the positivity of the dual variable suggests that a barrier method
may perform well. We leave this and other numerical exploration for future work.

\printbibliography

\ifsubmit\else
\appendix
\section{Simple examples}
In this section, we provide a number of simple examples using the \solver{} 
interface.

\subsection{Optimal power flow.}\label{sec:interface-opf}
First, we return to the optimal power flow example
from~\cite[\S3.2]{diamandis2024convexflows}. We wish to find a cost-minimizing
power generation plan that meets demand over a network of generators and
consumers connected by transmission lines. Given problem parameters demand
\jlinl{d}, line capacities \jlinl{ub}, and graph adjacency matrix \jlinl{Adj},
the entire optimal power flow problem may be defined and solved in less than ten
lines of code:

\begin{jllisting}
# Parameters: demand d, graph Adj, upper bounds ub
obj = NonpositiveQuadratic(d)
h(w) = 3w - 16.0*(log1pexp(0.25 * w) - log(2))

lines = Edge[]
for i in 1:n, j in i+1:n
    Adj[i, j] ≤ 0 && continue
    push!(lines, Edge((i, j); h=h, ub=ub[i]))
end

prob = problem(obj=obj, edges=lines)
result = solve!(prob)
\end{jllisting}
In this example, we used the special function \jlinl{log1pexp} from the 
\texttt{LogExpFunctions} package, which is a numerically well-behaved 
implementation of the function $x \mapsto \log(1 + e^x)$. Since $h$ may be
specified as native Julia code, using non-standard functions does not introduce
any additional complexity.

In this case, the arbitrage problem has a closed-form
solution, easily derived from the first-order optimality conditions. With a
small modification, we can give this closed-form solution to the solver directly.
We write this closed form solution as
\begin{jllisting}
# Closed for solution to the arbitrage problem, i.e. the wstar that solves
#  h'(wstar) == ratio
function wstar(ratio, b)
    if ratio ≥ 1.0 
        return 0.0
    else
        return min(4.0 * log((3.0 - ratio)/(1.0 + ratio)), b)
    end
end
\end{jllisting}
We only need to modify the line in which we define the edges, changing it
to
\begin{jllisting}
push!(lines, Edge((i, j); h=h, ub=ub_i, wstar = @inline w -> wstar(w, ub_i)))
\end{jllisting}
After these modifications, we lose little computational efficiency compared to
specifying the solution to the arbitrage problem directly. However, if we
specify the problem directly, we may pre-compute the zero-flow region, inside
which a line is not used at all. Compare the code for this example with the code
for the fully-specified example in~\cite[\S6.1]{diamandis2024convexflows}, which
can be found at
\[
    \texttt{\url{https://github.com/tjdiamandis/ConvexFlows.jl/tree/main/paper/opf}}.
\]

\subsection{Trading with constant function market makers}
Similarly, we can easily specify the problem of finding an optimal trade given a
network of decentralized exchanges. Here, we assume that the constant function
market makers are governed by the trading function
(see~\cite[\S3.5]{diamandis2024convexflows} for additional discussion)
\[
    \phi(R) = \sqrt{R_1 R_2}.
\]
This trading function results in the gain function
\[
    h(w) = \frac{w R_2}{R_1 + w},
\]
which one can easily verify is strictly concave and increasing for $w \ge 0$.
Given the adjacency matrix \jlinl{Adj} and constant function market maker 
reserves \jlinl{Rs}, we may specify the problem of finding the optimal trade as
\begin{jllisting}
obj = Linear(ones(n));
h(w, R1, R2) = R2*w/(R1 + w)

cfmms = Edge[]
for (i, inds) in enumerate(edge_inds)
    i1, i2 = inds
    push!(cfmms, Edge((i1, i2); h=w->h(w, Rs[i][1], Rs[i][2]), ub=1e6))
    push!(cfmms, Edge((i2, i1); h=w->h(w, Rs[i][2], Rs[i][1]), ub=1e6))
end

prob = problem(obj=obj, edges=cfmms)
result = solve!(prob)
\end{jllisting}
Here, our objective function $U(y) = \ones^Ty$ indicates that we 
value all tokens equally.

\subsection{Market clearing}
Finally, we revisit the market clearing example
from~\cite[\S3.4]{diamandis2024convexflows}. In this example, the objective
function includes a constraint, so we must specify it directly. However, we may
still specify edge gain functions instead of specifying the arbitrage problems
directly.

Recall that the objective function is given by
\[
    U(y) = \sum_{i=1}^{n_b} c_i \log y_i - I(y_{n_b+1 : n_b + n_g} \ge -1),
\]
where $n_b$ is the number of buyers, $n_g$ is the number of goods, and $c \in
\reals^{n_b}_+$ is a vector of budgets. We will define a struct 
\jlinl{MarketClearingObjective} to hold the problem parameters and then define 
the methods \jlinl{U} to evaluate the objective, \jlinl{Ubar} to evaluate the 
associated subproblem~\eqref{eq:subproblem-U}, and \jlinl{∇Ubar!} to evaluate
the gradient of the subproblem. The full implementation of the first subproblem 
is below.

\begin{jllisting}
const CF = ConvexFlows
struct MarketClearingObjective{T} <: Objective
    budget::Vector{T}
    nb::Int
    ng::Int
    ϵ::T
end
function MarketClearingObjective(budget::Vector{T}, nb::Int, ng::Int; tol=1e-8) where T
    @assert length(budget) == nb
    return MarketClearingObjective{T}(budget, nb, ng, tol)
end
Base.length(obj::MarketClearingObjective) = obj.nb + obj.ng

function CF.U(obj::MarketClearingObjective{T}, y) where T
    any(y[obj.nb+1] .< -1) && return -Inf
    return sum(obj.budget .* log.(y[1:obj.nb])) - obj.ϵ/2*sum(abs2, y[obj.nb+1:end])
end

function CF.Ubar(obj::MarketClearingObjective{T}, ν) where T
    return sum(log.(obj.budget ./ ν[1:obj.nb]) .- 1) + sum(ν[obj.nb+1:end])
end

function CF.∇Ubar!(g, obj::MarketClearingObjective{T}, ν) where T
    g[1:obj.nb] .= -obj.budget ./ ν[1:obj.nb]
    g[obj.nb+1:end] .= 1.0
    return nothing
end
\end{jllisting}
We specify the utility that buyer $b$ gets from good $g$ as
\[
    h(w) = \sqrt{b + gw} - \sqrt{b}.  
\]
With the objective defined, we may easily specify and solve this problem as
before:
\begin{jllisting}
obj = MarketClearingObjective(budgets, nb, ng)
u(x, b, g) = sqrt(b + g*x) - sqrt(b)

edges = Edge[]
for b in 1:nb, g in 1:ng
    # ub arbitrary since 1 unit per good enforced in objective
    push!(edges, Edge((nb + g, b); h=x->u(x, b, g), ub=1e3))
end

prob = problem(obj=obj, edges=edges)
result = solve!(prob)
\end{jllisting}
See the documentation for additional details and commentary.

\section{Multi-period power generation example}\label{app:multi-period-code}
We create a hour-by-hour power generation plan for an example network with three
nodes over a time period of $5$ days. The first two nodes are users who consume
power and have a sinusoidal demand with a period of $1$ day. These users may
generate power at a very high cost ($\gamma_i = 100$). The third node is a
generator, which may generate power at a low cost ($\gamma_i = 1$) and demands
no power for itself. The parameters are defined as follows:

\begin{jllisting}
# Problem parameters
n = 3
days = 5
T = 24*days
N = n*T

d_user = sin.((1:T) .* 2π ./ 24) .+ 1.5
c_user = 100.0
d_gen = 0.0*ones(T)
c_gen = 1.0

d = vec(vcat(d_user', d_user', d_gen'))
c = repeat([c_user, c_user, c_gen], T)
obj = NonpositiveQuadratic(d; a=c)
\end{jllisting}

Next, we build a network between these three nodes. We create the transmission
line as in~\S\ref{sec:interface-opf}. Then we build the storage edges, which
`transmit' power from time $t$ to time $t+1$. We equip the second user with a
battery, which can store power between time periods with efficiency $\gamma =
1.0$. The network has a total of $360$ nodes and $359$ edges. 

\begin{jllisting}
# Network: two nodes, both connected to generator
function build_edges(n, T; bat_node)
    net_edges = [(i,n) for i in 1:n-1]
    edges = Edge[]

    # Transmission line edges
    h(w) = 3w - 16.0*(log1pexp(0.25 * w) - log(2))
    function wstar(η, b)
        η ≥ 1.0 && return 0.0
        return min(4.0 * log((3.0 - η)/(1.0 + η)), b)
    end

    for (i,j) in net_edges
        bi = 4.0
        for t in 1:T
            it = i + (t-1)*n
            jt = j + (t-1)*n
            push!(edges, Edge((it, jt); h=h, ub=bi, wstar=η -> wstar(η, bi)))
            push!(edges, Edge((jt, it); h=h, ub=bi, wstar=η -> wstar(η, bi)))
        end
    end

    # Storage edges
    ϵ = 1e-2
    wstar_storage(η, γ, b) = η ≥ γ ? 0.0 : min(1/ϵ*(γ - η), b)
    
    # only node 2 has storage
    for t in 1:T-1
        it = bat_node + (t-1)*n
        it_next = bat_node + t*n
        γi = 1.0
        storage_capacity = 10.0
        push!(edges, Edge(
            (it, it_next); 
            h= w -> γi*w - ϵ/2*w^2, 
            ub=storage_capacity, 
            wstar = η -> wstar_storage(η, γi, storage_capacity)
        ))
    end
    return edges
end
\end{jllisting}

With the hard work of defining the network completed, we can construct and solve
the problem as before. We solve this problem with BFGS, as L-BFGS does not
exhibit good convergence on our problem, which is consistent to the results
in~\cite{asl2021behavior}. The (almost) linear edges mean this problem is 
(almost) nonsmooth.
\begin{jllisting}
edges = build_edges(n, T, bat_node=2)
prob = problem(obj=obj, edges=edges)
result_bfgs = solve!(prob; method=:bfgs)
\end{jllisting}
\fi
\end{document}